\documentclass[12pt]{article}
\usepackage{amssymb, url, amsthm, amsmath}
\usepackage{amscd}

\usepackage{textcomp}
\def\boxit#1{%
{\hbox{\lower3pt\hbox{\vrule\vbox{\hrule\kern2pt%
\hbox{\kern2pt$#1$\kern2pt}\kern2pt\hrule}\vrule}}}}
\bigskip
\bigskip
\def\be{\begin{equation}}
\def\ee{\end{equation}}

\def\ve{\varepsilon}

\def\R{{\sf I\kern-.2em R}}
\def\N{{\sf I\kern-.2em N}}
\def\C{\kern.1em{\raise.47ex\hbox{$\scriptscriptstyle
$}}\kern-.40em{\sf C}}
\def\Z{{\sf Z\kern-.32em Z}}

\def\hat{\widehat}

\def\hat{\widehat}
\def\be{\begin{equation}}
\def\ee{\end{equation}}

\def\ve{\varepsilon}

\newtheorem{theorem}{\noindent Theorem}

\newtheorem{definition}{\noindent Definition}

\newtheorem{statement}{\noindent Proposition}

\urldef\vershik\url{vershik@pdmi.ras.ru}
\urldef\gorbulsky\url{gorbulsky@consultant.com}
\author {A.~M.~Vershik\thanks{%
St.~Petersburg Department of Steklov Institute of Mathematics. E-mail: \vershik. Partially
supported by the grants RFBR 05-01-00089, and NSh-4329-2006-1.} \and A.~D.~Gorbulsky\thanks{E-mail:
\gorbulsky.}}

\date{25.05.07}
\title{Scaled entropy of filtrations of $\sigma$-fields}

\begin{document}

 \maketitle

\begin{abstract}
We study the notion of the scaled entropy of a filtration of
$\sigma$-fields (= decreasing sequence of
$\sigma$-fields) introduced in \cite{V4}. We suggest a method for
computing this entropy for the sequence of $\sigma$-fields of pasts
of a Markov process determined by a random walk over the trajectories of
a Bernoulli action of a commutative or nilpotent countable group
(Theorems~5,~6). Since the scaled entropy is a metric invariant
of the filtration, it follows that the sequences of $\sigma$-fields of
pasts of random walks over the trajectories of Bernoulli actions of lattices
(groups ${\Bbb Z}^d$) are metrically nonisomorphic for different dimensions
$d$, and for the same $d$ but different values of the entropy of
the Bernoulli scheme. We give a brief survey of the metric theory of
filtrations, in particular, formulate the standardness criterion
and describe its connections with the scaled entropy and the notion of a tower of measures.
\end{abstract}

 \tableofcontents

\section{Introduction: filtrations of
$\sigma$-fields; standardness; classification}

We begin with recalling some general definitions.
A {\it Lebesgue}, or {\it Lebesgue--Rokhlin}, {\it space}
$(X,\mu)$ is a space with a probability measure
$\mu$ that is metrically isomorphic
($\bmod 0$) to the union of the interval $[0,\lambda)$, $\lambda \leq 1$,
with the Lebesgue measure and, possibly, at most countably many atoms
of positive measures that sum to $1-\lambda$.
We will be interested mainly in Lebesgue spaces with continuous measures.
A {\it measurable partition} of a Lebesgue space
$(X,\mu)$ is the partition of $X$ into the preimages of points
under a measurable map; without loss of generality we may assume that
this measurable map is a real-valued measurable function
$f:X \to \Bbb R$. A class of $\bmod0$ coinciding functions
determines a class of $\bmod0$ coinciding partitions; in what
follows, speaking about partitions,
we always mean these classes rather than individual partitions.
Recall that, by Rokhlin's theorem
\cite{R1}, with every measurable partition
$\xi=\{C_{\alpha}\}$ with elements
$C_{\alpha}$, $\alpha \in \cal A$,
we can associate a canonical system of measures, namely, the system
of conditional measures
$\{\mu^C\}$ on the elements
$C_{\alpha}=C$; the conditional measures are well defined
for almost all elements of
$\xi$, so that the canonical system of measures is well defined
$\bmod 0$. The metric classification of
$\bmod 0$ classes of measurable partitions
in terms of systems of conditional measures
is due to V.~A.~Rokhlin \cite{R2}.

A measurable partition $\xi$ determines, and is determined by, a $\sigma$-subfield
${\frak A}_{\xi}$ of the $\sigma$-field
${\frak A}(X, \mu)$ of all classes of measurable sets of the space
$(X,\mu)$, namely, the $\sigma$-subfield generated by the Lebesgue
sets of the corresponding measurable function.
The language of $\sigma$-subfields of
${\frak A}(X, \mu)$, traditionally used in the theory of random processes,
is equivalent to the more geometric language of measurable partitions,
which we will mainly use in what follows. The correctness of definitions
with respect to considering classes of
$\bmod 0$ coinciding objects is usually easy to check (see, e.g.,
\cite{R1, V1}).

On the set ${\cal P}(X)$
of classes of measurable partitions ($\sigma$-fields)
there is a natural partial ordering. In terms of $\sigma$-fields,
it is the ordering by inclusion, with respect to which
${\cal P}(X)$ is a lattice.\footnote{In terms of measurable
partitions, ``greater'' in the sense of this ordering means ``finer,''
so that the greatest partition is the partition (denoted by
$\varepsilon\bmod 0$) into separate points; and the trivial partition,
denoted by $\nu \bmod 0$, whose two elements are the empty set
and the whole space, is the smallest element of the lattice of partitions.
This ordering is opposite to that accepted in combinatorics, where the
greatest element of the lattice is the trivial partition.}
We study infinite decreasing
sequences of measurable partitions (or infinite decreasing sequences of
$\sigma$-fields). {\it In this paper, the term ``filtration''
is a synonym of the term ``infinite decreasing sequence of measurable
partitions'' or ``infinite decreasing sequence
of $\sigma$-fields.''} A filtration
$\Xi=\{\xi_n,\, n \in \Bbb N\}$ is called
{\it ergodic} if the intersection $\bigcap\limits_n \xi_n$
of its components is the trivial partition
$\nu$, i.e., $\bigcap\limits_n \xi_n=\lim\limits_{n\to \infty}
\xi_n=\nu$.\footnote{The intersection of $\sigma$-fields is defined literally,
but in the language of measurable partitions, the intersection is the
measurable hull of the individual (set-theoretic) intersection
of partitions.} A general example of a filtration is the sequence
of $\sigma$-fields of ``pasts'' of a one-sided discrete-time
random process $\{y_n , \, n\le0\}$, i.e., the sequence
$\{{\cal A}_n\}_{n=0}^{\infty}$, where ${\cal A}_n$ is
the $\sigma$-field generated by the random variables
$\{y_k : k\le -n\}$. This filtration is ergodic if and only if the infinite
past is trivial (i.e., the process is Kolmogorov-regular). It is of special interest
to study the
{\it sequences of pasts of stationary random processes} considered below;
in this case, the sequence of partitions is shift-invariant (or, in short,
stationary). For more details on this theory, see
\cite{V1} and the references therein. Stationary filtrations (i.e., the
sequences of pasts of stationary discrete-time or continuous-time processes)
is one of the two main objects
of filtration theory. The second class
of examples, which is not less important, consists of filtrations
arising in trajectory theory and the theory of periodic approximations of
dynamical systems; here we do not consider this class. From the
point of view of the theory of stationary random processes,
the filtration of pasts, its structure and its metric type, is the most
important characteristic of the process and contains
deep information about it.

Filtrations $\Xi=\{\xi_n\}_{n=1}^{\infty}$ and
$\Xi'=\{\xi'_n\}_{n=1}^{\infty}$ are called (metrically) {\it isomorphic}
if there exists a measure-preserving measurable transformation $T$
satisfying the condition $T\xi_n = \xi'_n$ for all $n$. The problem
of metric classification of infinite ergodic filtrations
was posed by the first author (mainly in connection with trajectory
theory) and has accumulated much literature.

The simplest example of a filtration is the Bernoulli filtration
which consists of the pasts
of a one-sided stationary Bernoulli scheme. It is ergodic, as follows from Kolmogorov's
zero or one law. A filtration metrically isomorphic to a Bernoulli
filtration is called {\it standard}; it is determined
by the type of the one-dimensional distribution of the Bernoulli scheme.
The Bernoulli scheme with probabilities
$1/r, \ldots, 1/r$, $r>2$, $r
\in \Bbb N$, determines a {\it standard $r$-adic filtration};
if $r=2$, a {\it dyadic filtration};
if the one-dimensional distribution of the Bernoulli scheme
is continuous,
a {\it standard continuous filtration}. More general
{\it nonstationary $\{r_n\}$-adic standard filtrations}
arise from nonstationary Bernoulli schemes. Filtrations
$\Xi=\{\xi_n\}$ and
$\Xi'=\{\xi'_n\}$ are {\it finitely isomorphic} if their finite fragments
$\{\xi_n\}_{n=1}^m$ and
$\{\xi'_n\}_{n=1}^m$
are isomorphic for any length $m$.
A filtration $\{\xi_n\}$ that is finitely isomorphic to
a standard $r$-adic (respectively, dyadic, continuous,
$\{r_n\}$-adic) filtration
is called {\it homogeneous $r$-adic}
(respectively, {\it dyadic, continuous,
$\{r_n\}$-adic});  and a general {\it homogeneous} filtration is a
filtration that is finitely isomorphic to an arbitrary (possibly,
nonstationary) standard filtration.

The original question was {\it whether finitely isomorphic
homogeneous ergodic filtrations can be nonisomorphic};
in other words, whether they can be essentially different ``at infinity''
provided that all their finite fragments are isomorphic?
For example, whether there exist metrically nonisomorphic ergodic dyadic
filtrations? The positive answer to this question, and thus the first example of a
nonstandard ergodic dyadic sequence, was obtained in
\cite{V0} (the detailed proofs were presented in
\cite{V2, V1}); this example is
the sequence of pasts of a random walk over
the generators of a Bernoulli action of the free group with
two generators. This example and its further generalizations
showed, in particular, that the metric type of the filtration
of pasts of a stationary process can be essentially different for
different stationary processes, and the corresponding classification problem
is meaningful. The first method for distinguishing filtrations
was combinatorial, but in fact it was of entropy nature. It led
to the definition of the combinatorial (or exponential) entropy
of a filtration (see \cite{V3} and below). This made it possible to
present a continuum of pairwise nonisomorphic dyadic filtrations.
In \cite{S} it was observed that the entropy of the action of the
dyadic group $\sum {\Bbb Z}_2$ is also an invariant of the filtration
generated by this action, and this also gives a continuum of nonisomorphic
filtrations. Moreover, in the dyadic case, the combinatorial entropy
and the entropy of the action coincide, though their definitions are
quite different. The coincidence of these entropies even for
$\{r_n\}$-adic sequences holds only for a certain growth of the number
$\{r_n\}$ of points in the elements of the partitions (see
\cite{V1, He1}). Besides, entropy of action can be defined only
for homogeneous filtrations, while combinatorial entropy is defined
for general filtrations (see below).

Combinatorial entropy distinguishes only a very narrow class
of filtrations, namely, filtrations with exponential asymptotics of the
iterated semimetrics (see below). Later, in
\cite{V4}, the class of {\it scaled entropies},
which we deal with in this paper, was introduced as
a generalization of
the notion of combinatorial entropy.
The definition of scaled entropy is based on introducing a scaling
for the growth of the entropies of appropriate partitions. In the hierarchy
of these scalings, combinatorial entropy corresponds to
exponential growth, so that it can be called
exponential entropy.

All currently known results illustrate the fact, unobvious {\it a priori},
that the metric classification of general, or even stationary, filtrations
is as difficult
as, e.g., the metric classification of stationary processes themselves.
This is exactly why the problem of finding constructive metric invariants
of homogeneous ergodic filtrations arises.
Entropy techniques, discussed in this paper, seem most useful in this regard.
Among other general theorems on filtrations, we would like to mention
the theorem on lacunary isomorphism and the ensuing notion of the
principal invariant of
filtrations (see \cite{V1});
the standardness criterion
suggested in \cite{V0, V2} for distinguishing between standard and nonstandard
filtrations (see \S\,2) is also partly motivated by
this theorem. {\it Scaled entropy}
is exactly the quantitative characteristic of filtrations that naturally
arises from the analysis of this criterion. It is determined by a so-called
scaling function (see below). In this paper we formulate
theorems on the scaled entropy of the filtrations of pasts of random
walks over the trajectories of Bernoulli actions of countable commutative or
nilpotent groups and outline their proofs. Possibly, this method applies
to groups for which the central limit theorem for random walks holds.
For an abelian or nilpotent group, the scaling is the power function
$n^{d/2}$, where $d$ is the weighted dimension of the group; in particular,
for the lattice ${\Bbb Z}^d$ it is equal to
$n^{d/2}$ (Theorems 5, 6). Moreover, it turns out that not only the order
(scaling), but also the value of the scaled entropy is an invariant. Thus
the filtrations of pasts of random walks on the lattices
${\Bbb Z}^d$ are metrically nonisomorphic for different dimensions $d$,
and even for the same $d$ but different values of the average entropy
of the transition probabilities.

The analysis of the filtrations of pasts of stationary processes provides
new characteristics of one-sided processes. As we will see, already for
Kolmogorov-regular processes, i.e., processes with trivial infinite past,
the filtrations of pasts can have quite various metric properties.
It is also possible that some invariants of the filtration of pasts
of a random process can be invariants of the two-sided shift in the space
of trajectories of the process. Considering the scaled entropy of filtrations
arising in problems of periodic approximation of automorphisms
leads to new invariants, such as
the scale of an automorphism and the so-called principal invariant,
see \cite{V5}.\footnote{Note that
the notion, introduced in \cite{V4}, of the {\it secondary entropy} of a stationary random process
is close to the notion
of scaled entropy; a similar characteristic was also studied in
\cite{Sh}.}

\section{Iterated Kantorovich metric, standardness criterion, and
the tower of measures}

\subsection{Admissible metrics and the Kantorovich distance}

In order to construct invariants of filtrations and, in particular,
formulate the standardness criterion, we need the
construction of {\it iterated metrics} and the notion of {\it tower of
measures}. But first we give the definitions of admissible metrics on
a measure space (admissible triples) and the classical Kantorovich
metric on measures.

\begin{definition}
We say that a semimetric $\rho$ on a Lebesgue space $(X,\mu)$ is {\em
admissible} (or the triple $(X,\mu,\rho)$ is admissible)
if the following conditions hold:

{\rm1)} the semimetric
$\rho(x,y)$, regarded as a function of two variables (i.e., as
a function on the space $(X\times X,\mu \times \mu)$), is measurable;

{\rm 2)} in the space $X$ there exists a subset $X_0$
of full measure $\mu$ that is quasi-compact, in the topological sense,
with respect to $\rho$; this means that the quotient space
${\hat X}_0$ of $X_0$ with respect to the partition into classes
of points with pairwise zero distances, endowed with the quotient metric,
is a compact metric space.
\end{definition}

As above, we consider classes of metrics (semimetrics)
coinciding almost everywhere rather than individual metrics (semimetrics).
Denote the set of all (classes of) admissible metrics on a given
Lebesgue space $(X,\mu)$ by $\Psi(X,\mu)$. The metric compact triple
$(X,\rho, \mu)$, where $\rho$ is a metric that turns
$X$ into a compact metric space
and $\mu$ is a probability Borel measure on this space,
is an example of an admissible triple.

Now recall the definition of the Kantorovich metric on the space of
measures on a compact metric space.

Given a compact metric space $(X, \rho)$,
one can define the Kantorovich metric
$k_{\rho}$ on the simplex $V(X)$ of probability Borel measures on $X$ (see
\cite{K}, and also a modern exposition \cite{V7}). The classical
definition of the Kantorovich metric applies only to compact metric spaces,
but it can be extended, without essential changes,
to the case of semimetrics and quasi-compact spaces.
This definition is as follows:
   $$k_{\rho}(\mu_1,\mu_2)=\inf\{\int\limits_{X\times X} \rho(x,y)\, d
Q(x,y)\mid P_1Q=\mu_1, P_2Q=\mu_2 \};$$
here $Q$ ranges over the set of all probability measures on
$X \times X$ with the given projections, $\mu_1$ and $\mu_2$, to both coordinates,
or, in the accepted terminology, with the given
marginal distributions; and
$P_1$ and $P_2$ are the projections which map measures on $X\times X$
to the simplices of measures on the corresponding coordinates.

\subsection{Iterated semimetrics associated with a filtration, and
the standardness criterion}

Now let us apply Kantorovich's construction to
constructing metrics on a measure space with
a given measurable partition. The following procedure,
suggested in \cite{V1, V2}, allows one,
{\it given an admissible semimetric $\rho$ and a measurable partition
$\xi$, construct a new admissible semimetric on the same space
$(X, \mu)$}. Let us fix an admissible semimetric
$\rho=\rho_0$ on the space $(X,\mu)$ and define a new distance
$\rho_1(x,y)$ on $(X, \mu)$
as the Kantorovich distance between the conditional measures
on the elements of $\xi$ that contain the given points:
$$ \rho_1(x,y)\doteq k_{\rho}(\mu^{C(x)}, \mu^{C(y)}),$$
where $C(x),C(y)$ are the elements of $\xi$ that contain
$x$ and $y$, respectively, and
$\mu^C$ is the conditional measure on an element
$C \in \xi$.

{\it Thus, given a semimetric and a measurable partition, we can
define a new semimetric.} The new distance between points
lying in the same element of the partition is equal to zero,
so that $\rho_1$ is a semimetric even in the case when
$\rho$ is a metric. However, the quotient of this semimetric on the
quotient space $X/\xi$ is a well-defined metric.

Now assume that we have a space $(X,\mu)$ with an admissible semimetric
$\rho$ and a filtration $\Xi=\{\xi_n\}_1^{\infty}$, $\xi_0=\varepsilon$.
Let us successively apply the above procedure:
using $\rho$ and $\xi_1$, construct a semimetric $\rho_1$;
then, using $\rho_1$ and $\xi_2$, construct a semimetric $\rho_2$, etc.
Since the partitions decrease, the semimetrics are coherent, in
the sense that if the distance between two points vanishes with
respect to
$\rho_k$, then it vanishes with respect to all subsequent semimetrics.

Note that the average distance between pairs of points of the space
$X$ does not increase when passing to the next iteration, i.e.,
$$
c_n=\int\limits_{X\times X}\rho_n(x,y)\le\, c_{n-1}=\int\limits_{X\times X}\rho_{n-1}(x,y)d\mu(x)
d\mu(y).
$$

Now we can formulate the standardness criterion for a homogeneous filtration.

\begin{theorem}[Standardness criterion, \cite{V0,V1,V2}]
A homogeneous filtration (a homogeneous sequence of partitions)
is standard if and only if for every initial semimetric
$\rho$ the mean value of the iterated distance between points tends to zero.
In our notation, the latter condition reads as
$$
\lim_{n\to \infty}\int\limits_{X\times X}\rho_n(x,y)d\mu(x)d\mu(y)=0, \quad
\mbox{i.e.,}\quad
\lim_n c_n =0.
$$
\end{theorem}

By definition, the condition of this criterion is metrically invariant.
Hence the ``only if'' part immediately follows from the
fact that it is satisfied for a Bernoulli filtration.
The ``if'' part is a deep result. The standardness
criterion was formulated in \cite{V0} and proved in full detail for
dyadic filtrations in \cite{V2}; the proof was reproduced
for $r$-adic filtrations  in
\cite{V1}.

For $r$-adic filtrations, the standardness criterion has a clear
combinatorial meaning. In this case, it suffices to check it
for semimetrics that reduce
to finite metric spaces. If we start
with such a semimetric $\rho$ (let it reduce to a $k$-point space),
then the $n$th iterated semimetric is a semimetric on the orbits of the action of
the group $D_{n,r}$ of automorphisms of the homogeneous one-root tree
$T_{n,r}$ of height $n$ and valence $r$
on the space of functions on
$T_{n,r}$ with values in a $k$-point set (for instance,
$\textbf{k}=\{0,1, \dots, k-1\}$). In more detail, the group $D_{n,r}$
acts by substitutions on the cube $\textbf{k}^{r^n}$ endowed
with the Hamming metric, and the iterated metric
$\rho_n$ reduces to the space of orbits of this group endowed
with the quotient metric (the distance between two orbits is
the minimum distance between the points of these orbits).
It is difficult to compute  this metric explicitly, but
in many cases it is possible to check whether or not it degenerates
in the limit (i.e., whether or not the space reduces to the one-point space).
For example, this can be done for
$r=k=2$. It is this computation that led to the first example of a nonstandard
filtration. It was carried out in
\cite{V1,V2} for the random walk over the trajectories of a Bernoulli
action of the free group and for a Bernoulli action of the group of
$2$-adic integers.

To prove the criterion
for filtrations with continuous conditional measures, one may use
the same scheme as in \cite{V1},
making only minimal changes compared with
the case of dyadic or $r$-adic filtrations.
Later, other proofs were suggested
for the continuous case, which have essentially the same ideology
as in the discrete case, see
\cite{DFST}. A detailed survey of questions related to standardness
and other properties of filtrations considered by B.~Tsirelson is given in
\cite{Em}; the latter paper also contains a rather complete list of
references.

Below we refine the standardness criterion. Namely, we want not only to know
whether or not the metric degenerates, but also to obtain an estimate on
the asymptotics of the entropy of the compact metric space
with the iterated metric. This
is the next step in the study of nonstandard filtrations.

\subsection{Tower of measures}

\def\limind{\operatorname{lim\,ind}\limits}
\def\limproj{\operatorname{lim\,proj}\limits}
\def\Tow{{\rm Tow}}

Another formulation of the standardness criterion uses the concept,
important in itself, of {\it tower of measures}
\cite{V1}, which we will briefly reproduce here (see also
\cite{V3}). First assume that we are given a compact metric space
$(Y,r)$; consider the simplex $V(Y)$ of probability Borel measures on $(Y,r)$
endowed with the Kantorovich metric
$k_r$; it is also compact in the topology determined by this metric
(i.e., in the weak topology). Moreover, there exists an isometric
embedding $i:y\to \delta_y$ of the initial space $Y$ into
$V(Y)$. Then we can consider the isometric embedding of the simplex
$V(Y)$ into the simplex $V(V(Y))\equiv V^2(Y)$
of probability Borel measures on $V(Y)$, again with the Kantorovich metric
$k_{k_r}$, and so on. Thus we have an inductive family of compact
metric spaces $V^n=V^n(Y)$
with isometric embeddings $i_n: V^n\to V^{n+1}$,
and we can consider the inductive
limit of these spaces:
$$
\limind_n (V^n(Y),i_n)=V^{\infty}(Y)\equiv \Tow(Y,r);
$$
the limit space (inductive limit) $\Tow(Y,r)$ is called the
{\it tower of measures}; it is endowed with an inductively
defined metric, which we
denote by
$\bar r$. This definition is of purely topological nature;
the construction can be generalized to the case of a
quasi-compact semimetric space, but we will not need such a generalization. The space
$(\Tow(Y,r), \bar r)$ is a (noncomplete) metric space. Its nature is
worth a detailed study; it is of special interest and importance to study
the properties of its completion with respect to
the metric $\bar r$.
Note that it has not only the structure of an inductive limit, but
also the structure of a projective limit. Indeed, since
$V^n(Y)$, $n>1$, is an affine compact space, every measure from
$V^n(Y)$ has a well-defined barycenter, which is
a point of $V^{n-1}(Y)$. The map that sends a measure
to its barycenter is an epimorphic affine projection
$V^n(Y) \to V^{n-1}(Y)$, $n>1$, right inverse to the embedding
$i_{n-1}$. This allows one to consider the projective
limit of compact spaces $\limproj V^n$; the limit compact space is exactly
the completion of the inductive limit.\footnote{It is natural to
say that a space with
such coherent structures of an inductive limit and a projective limit
has the structure of an ``indoprojective limit.''}

\subsection{The standardness criterion in terms of the tower of measures}

Let us apply the tower of measures construction to the study of filtrations. Assume
that we are given a Lebesgue space
$(X,\mu)$ and an arbitrary measurable function
$f:X \to [0,1]$. Set
$f_0\equiv f$ and consider the tower of measures
$\Tow([0,1], r)$, where $r$ is the Euclidean metric on
$[0,1]$. Assume that we are given a filtration
$\{\xi_n\}_{n=0}^{\infty}$ on $(X,\mu)$. Let us define a sequence of
probability measures $\{\nu_f^n\}$, $n=1,2, \dots$, on
$\Tow([0,1], r)$, where
$\nu_f^n\equiv\nu^n$ is a measure on
$V^n \subset \Tow([0,1], r)$, as follows.
The first measure
$\nu^1 \in V^1([0,1])\subset \Tow([0,1], r)$ is the image of $\mu$
under the map $f_0:X \to [0,1]$; thus
$\nu^1$ is a measure on $[0,1]$, i.e., an element of $V^1([0,1])$.
Then we consider the map
$f_1: x \mapsto f_0(\mu^{C_1(x)})\in V^1$ that sends a point
$x\in X$ to the image of the conditional measure on the element
$C_1(x)$ under the map $f_0$ restricted to $C_1(x)$.
The second measure $\nu^2$ is the image of
$\mu$ under $f_1$, i.e., a measure on $V^1([0,1])$
(a ``measure on measures'' on $[0,1]$, or an element of
$V^2([0,1])$). Note that the function
$f_1$ is well defined on the quotient space
$X/\xi_1$. Now we consider the map
$f_2:X/\xi_1 \to V^2([0,1])$ that sends a point $y\in X/\xi_1$ to
the image of the conditional measure on the element of
$\xi_2/\xi_1$  containing $y$
under the function $f_1$ restricted to this element.
The measure $\nu^3$ is the image of $\mu/\xi_1$
under $f_2$, so that it is a measure on
$V^2([0,1])$ (or an element of
$V^3([0,1])$), and so on. In this way we inductively define a map
$f_n:X\to X/\xi_{n-1}\to
V^n([0,1])$, $n=1,2,\dots$, and a measure $\nu^{n+1}$, which
is the image of $\mu$ under $f_n$.  Thus we have constructed a sequence of
measures $\nu^n$, $n=1,2, \dots$, on compact spaces lying in the tower of measures
$\Tow([0,1],r)$.  For more details, see
\cite{V1,Em}.\footnote{In \cite{V1}, the map $f_n:X \to
V^n([0,1])$, more exactly, the map that sends the initial function
$f_0=f:X \to [0,1]$ to $f_n$, was called the {\it
universal projection} of $f$ with respect to the finite decreasing
sequence of partitions
$\Xi_n=\{\xi_k,\,k=1,2,\dots, n\}$;
all joint metric invariants of $f$ and
$\Xi_n$ can be expressed in terms of the functions $f_n$.}

In these terms, the standardness criterion asserts that a filtration
is standard if and only if for every measurable function
$f=f_0$, the sequence of measures
$\nu_f^n$ collapses to a delta measure on the completion of the
tower of measures (i.e., the weak limit of
$\nu_f^n$ is the delta measure at a point belonging to the completion
of the tower of measures). This can be expressed by the following formula:
$$
\lim_{n\to\infty}\iint{\bar r}(x,y)d \nu^n(y) d \nu^n(y)=0,
$$
where the integral is taken over the square of
$\Tow([0,1], r)$ and ${\bar
r}$ is the metric on  $\Tow([0,1], r)$
defined by the above tower of measures construction applied to
the space $([0,1],r)$ with $r$ the Euclidean metric.
Compared to the first formulation of the standardness criterion (see above),
the integration of the iterated metric over the space
$(X,\mu)$ is replaced in this formula
by the integration over the tower of measures.

In fact, the standardness takes place if the above condition
is fulfilled for at least one one-to-one $\bmod 0$ measurable
function. The fact that a filtration is not standard means, on
the contrary, that there exists a function for which  the sequence of
measures $\nu^n$ does not degenerate, and this function is not measurable
with respect to any coherent sequence of independent complements to
the filtration
(see \cite{V2}). This interpretation immediately implies that the behavior
of the sequence of measures $\nu^n$ (or,
in the first interpretation, the sequence of metrics  $\rho_n$ on
$(X,\mu)$) contains essential information
on the asymptotics of the filtration. In fact, the metric type
of the filtration is determined by the sequence of measures
$\{\nu_f^n\}_{n=1}^{\infty}$ on $\Tow(X,\rho)$
associated with a one-to-one $\bmod0$ function $f$. In particular,
the asymptotics of the $\varepsilon$-entropy of the measures
$\nu^n$, $n=1,2, \dots $, on
$\Tow(X,\rho)$ is an invariant of the filtration and does not depend
on the choice of the initial metric $\rho$.

One can study filtrations either in terms
of the iterated metrics
$\rho_n$ (which will be done below) or in (equivalent) terms of the measures
$\nu_n$ on the tower of measures. It is the interrelation between the
filtration and the $\varepsilon$-entropy of the metric measure spaces
$(X,\mu,\rho_n)$
that will be used in the next section for introducing the notion of scaling
and scaled entropy. In brief, the difference between the two formulations
of the standardness criterion can be expressed as follows:
in the first case, we fix the measure space and iterate the
metric; in the second case, we fix the compact metric space and the
associated tower
of measures and vary measures on the tower of measures.
Apparently, these two approaches are equivalent not only in the
formulation of the standardness criterion, but also in the analysis
of numerical characteristics constructed from the metrics
$\rho_n$ in the first case and the measures
$\mu_n$ in the second case.

\section{Scaled entropy of filtrations: definition, examples}

\subsection{Entropy of a metric measure space}

The ordinary definitions of the
$\varepsilon$-entropy of a compact metric space and the entropy of an atomic measure
are well known (see, e.g.,
\cite{R1}). We will need the following characteristic
of a metric measure space.

\begin{definition}
Let $\ve>0$. The
$\varepsilon$-{\em entropy} of a semimetric measure space
$(X,\mu,\rho)$, with $\rho$ an admissible semimetric, is the following function of
$\varepsilon$:
$$H_{\varepsilon}(X,\rho,\mu)=\inf \{H(\lambda)\mid k_{\rho}(\lambda ,\mu) < \varepsilon\},$$
where $\lambda$ ranges over all discrete measures on
$X$, $H(\cdot)$ is the entropy of a discrete measure, and
$k_{\rho}$ is the Kantorovich metric on the space of measures on $X$.
\end{definition}

Roughly speaking,
$H_{\varepsilon}(X,\rho,\mu)$ is the ``entropy'' of the continuous
measure $\mu$ in the semimetric space up to
$\varepsilon$.\footnote{In the literal sense, the entropy of
a continuous measure is equal to infinity, and the definition
of what ``up to'' means depends on the semimetric $\rho$.}

We will use the $\ve$-entropy of a space
$(X,\mu)$ endowed with an ergodic filtration and the associated sequence of
admissible iterated semimetrics
$\rho_n$. The analysis of the asymptotic behavior of this
$\ve$-entropy allows us to define and compute the scaled entropy of the filtration.

\subsection{The definition of scaled entropy}

\begin{definition}
We say that a positive function $c(\varepsilon,n)$ of two arguments
$\varepsilon >0$, $n \in \Bbb N$ is a {\em scaling function}
if it is increasing in $n$ for a fixed
$\varepsilon$ and nonincreasing in $\ve$ for a fixed $n$.
Two scaling functions
$c(\cdot,\cdot)$ and $c'(\cdot,\cdot)$ are {\em strictly equivalent}
if $$
\overline{\lim\limits_{\varepsilon\to 0}}\quad
\overline{\lim\limits_{n\to \infty}} {c(\varepsilon,n)\over
c'(\varepsilon,n)} =\overline{\lim\limits_{\varepsilon\to 0}}\quad
\underline{\lim\limits_{n\to \infty}} {c(\varepsilon,n)\over
c'(\varepsilon,n)}=1.
$$
If each of these limits is equal to a finite nonzero number,
then the scaling functions $c(\cdot,\cdot)$ and $c'(\cdot,\cdot)$
are called {\em equivalent}.
\end{definition}

\begin{definition} The {\em scaled entropy} of a filtration
$\{\xi_n\}=\Xi$ with respect to a semimetric
$\rho$ with scaling function $c(\cdot,\cdot)$ is
the number
$$
h_c(\Xi,\rho)=\limsup\limits_{\varepsilon \to 0}\limsup\limits_{n\to \infty}
{H_{\varepsilon}(X,\mu,\rho_n)\over c(\varepsilon,n)},
$$
where $\rho_n$ is the iterated semimetric associated with the filtration
$\{\xi_n\}=\Xi$ (see the previous section). A {\em proper scaling function}
of $\Xi$ is a scaling function for which the scaled entropy
$h_c(\Xi,\rho)$ is different from zero and infinity.
\end{definition}

The existence of a proper scaling function is a separate problem. In the
examples below it is easy to prove.

The following proposition is obvious.

\begin{statement}
{\rm 1.} The values of the scaled entropy
$h_c(\Xi,\rho)$ with respect to a given semimetric $\rho$
with strictly equivalent scaling functions coincide.

{\rm 2.} For any filtration and a given semimetric there exists at most
one, up to equivalence, proper scaling function.
\end{statement}

Note that we may compute the scaled entropy
$h_c(\Xi,\rho)$ of a filtration with an arbitrary scaling function,
but the answer will be different from zero and infinity for at most
one class of equivalent scalings. Sometimes, for
filtrations of a certain type (for example, $r$-adic)
one can choose one distinguished normalization scaling
(similarly to choosing the base of logarithms in the definition
of the ordinary entropy).
In this case, the common value of the scaled entropy with respect to
a given semimetric for all proper scaling functions strictly equivalent
to the normalization scaling is called the value of the scaled entropy
with respect to the given semimetric.

The definition of scaled entropy is as follows
(see \cite{V4}).

\begin{definition}
Assume that for some class of filtrations we have chosen a normalization scaling. The {\em scaled
entropy} of a filtration $\Xi$ is the supremum of the normed scaled entropies with respect to
$\rho$ over all admissible semimetrics $\rho$, i.e., the following finite or infinite number:
$$
h(\Xi)=\sup_{\rho} h(\Xi,\rho).
$$
\end{definition}

Recall that $\rho$ is the initial semimetric from which the iterated
semimetrics
$\rho_n$ were constructed.

\begin{theorem}
The numerical value $h(\Xi)$ of the scaled entropy (if a normalization
scaling exists) is a metric invariant of the filtration
$\Xi$.
\end{theorem}

Although, as follows from this theorem, it suffices to start
the construction of iterated semimetrics from a metric, nevertheless we
consider semimetrics, since it is more convenient for calculations
and in this way it is easier to approximate the value
of the scaled entropy.

When calculating the scaled entropy of an arbitrary filtration,
it is natural to start searching for a correct scaling with the
exponential scaling (see below) and, if the corresponding entropy
vanishes, turn to another scaling with slower growth, and so on.
For a standard filtration, the entropy vanishes for any scaling.
The converse is also true, because, in view of the standardness criterion,
the standardness means that the scaled entropy vanishes for the scaling
$c(\ve,n) = {\rm const}$. This fact is similar to Kushnirenko's theorem
\cite{Ku} on the action of $\Bbb Z$, which states that
the vanishing of all sequential
entropies of an automorphism is equivalent to the discreteness of its spectrum.
The calculation of the scaled entropy of a filtration
is an interesting and not very simple problem. It is not known even what
scaling functions can really appear in the definition of the scaled entropy
of filtrations. An interesting question concerns the connections of the
scaled entropy of, e.g., dyadic sequences to properties
of group actions. Below we find the scaling and calculate the entropy
for the filtrations of pasts of several important Markov processes.

\subsection{Exponential (combinatorial) entropy}

The concept of scaled entropy arises, on the one hand,
from the analysis of the standardness
criterion and, on the other hand, from the original notion of the entropy
of a filtration suggested in \cite{V2}, which, from the viewpoint of
the above definition, is the scaled entropy with exponential
scaling. We call it the {\it exponential entropy}. Let us briefly
provide some information on this entropy;
for definiteness, we restrict ourselves to the case of homogeneous filtrations.

We will consider $\{r_n\}$-adic homogeneous filtrations (see Sec.~1).
Almost every element of the partition
$\xi_n$ of an $r_n$-adic filtration consists of
$\prod\limits_{i=1}^n r_i$ points; on each element the previous
partitions determine the structure of a (hierarchy) tree. Denote the
group of automorphisms of this tree by
$D_{\{r_k\}}$, $k=1,\dots, n$. If we fix an arbitrary finite partition
$\gamma$ and label its elements, in an arbitrary way, by some symbols
$0,1, \dots, p$, then for every $n$, for every element of $\xi_n$
we can define a sequence of length
$\prod\limits_{i=1}^n r_i$ whose coordinates are the symbols $0,1,
\dots, p$ corresponding to the points from the given element of
$\gamma$. The definition of this sequence
is not invariant (it depends on the labelling of points in the element of
the partition), but the orbit of the action of the group
$D_{\{r_k\}}$ on such sequences already depends only on the point
$x$ and the partition $\gamma$ (and, of course, on the fragment of length $n$
of the filtration). Thus we obtain a partition
$\gamma_n$ whose every element consists of all points having the
same orbit of the action of the group
$D_{\{r_k\}}$ on sequences consisting of the symbols
$0,1, \dots, p$.

\begin{definition}
The {\em entropy of an $r_n$-adic sequence of partitions
$\Xi=\{\xi_n\}_n$ with respect to a finite partition} $\gamma$
is the number
$$
h(\Xi;\gamma)=\lim\limits_{n\to \infty}{1\over \prod\limits_i r_i}
H(\gamma_n),
$$
where $H(\cdot)$ is the binary entropy of a finite partition.
\end{definition}

Note that  $H(\gamma_n)\le r_n\cdot H(\gamma_{n-1})$, so that
$h(\Xi;\gamma)$ is bounded.

Now we can get rid of
$\gamma$ and define an invariant of the filtration.

\begin{definition}
The {\em exponential entropy of an
$r_n$-adic sequence of partitions} is the number
$$
h(\Xi)=\sup\limits_{\gamma} h(\Xi;\gamma).
$$
\end{definition}

The constructed invariant can also be called the combinatorial entropy
of a homogeneous filtration. The scaling depends on the sequence
$\{r_n\}$.  The following theorem, which is an analog of
Kolmogorov's theorem, was proved in \cite{V2}.

\begin{theorem}
The entropy $h(\Xi;\gamma)$ is continuous in
$\gamma$, with respect to the metric
$H(\gamma_1,\gamma_2)=H(\gamma_1|\gamma_2)+H(\gamma_2|\gamma_1)$
on the space of finite partitions.
\end{theorem}

This allows one to approximate the exponential entropy by
the values $h(\Xi;\gamma)$ for appropriate partitions $\gamma$.
An easy consequence of Theorem~3
is the following fact:
\textit{the exponential entropy of a standard filtration is equal to zero}.

For filtrations generated by actions of locally finite groups,
it is natural to compare the exponential entropy with the entropy
of the action (see \cite{S}); they coincide if the growth of
$\{r_n\}$ is not too fast (see \cite{V5, He1}).

The following theorem from
\cite{G2} shows that exponential entropy is a special case of
scaled entropy.

\begin{theorem}
The exponential entropy of an
$\{r_n\}$-adic filtration coincides with the scaled entropy
with the scaling function
$c(\varepsilon,n)=\prod\limits_{i=1}^n r_i$.
\end{theorem}

Exponential entropy vanishes for a wide class of filtrations
(see the next section), so that it does not solve
the problem of classification
of decreasing sequences of partitions. This is demonstrated by
a number of examples
of nonstandard filtrations with zero exponential entropy.
Scaled entropy allows one to further distinguish metrically
nonisomorphic filtrations.

If the scaled entropy does not vanish for some nonexponential scaling,
then the exponential entropy of such a filtration vanishes. In other words,
for an $\{r_n\}$-adic filtration,
the exponential scaling with the normalization chosen above
is maximal possible up to equivalence.

Exponential entropy can be defined not only for
$\{r_n\}$-adic filtrations,
but for arbitrary filtrations, including those
with continuous conditional measures. We do not dwell
upon this question.

It is well known that in the ordinary entropy theory, introducing
a scaling for measuring the rate of growth (in $n$) of the entropy of the
product of $n$ shifts of a partition in the case when the entropy
of the shift vanishes, does not lead to new invariants.
The reason is that for every ergodic automorphism with zero entropy,
one can choose the initial partition so that the growth of the entropies
of the product of
$n$ shifts of this partition as $n$ tends to infinity will have
a given subexponential rate. In the definition of the scaled
entropy of a filtration, we avoid this difficulty by
using admissible metrics. This allows us to construct invariants for
distinguishing
different asymptotics of the growth of the entropies.
For a given metric, our definition could also be compared
with the definition of the topological entropy of transformations,
but the significant difference is that we then take the supremum over
admissible metrics. Thus one can conjecture that the idea of
scaled entropy can be used also for automorphisms with zero
entropy.\footnote{Usually one fixes the metric and varies the measure
(cf.\ the notion of measure of maximal entropy); but we, on the
contrary, fix the measure and vary the metric. This idea
was repeatedly used in the papers of the first author.}

\section{The scaled entropy of filtrations generated by random walks
over the trajectories of group actions}

\subsection{Standardness and random walks}

As mentioned above, the first example of a nonstandard filtration
was the dyadic sequence of pasts of the random walk over the trajectories
of a Bernoulli action of the free group with two generators \cite{V1}.
The proof consisted in calculating the lengths of orbits
of the group of automorphisms of the hierarchy on the elements of the partitions.
Namely, let $\nu^n$ be the measure on the tower of measures
constructed from the semimetric corresponding to the
characteristic function of some set.
It was proved that there is no orbit of such a length that
the measure $\nu^n$ is concentrated in a neighborhood of this orbit.
Merely it was shown that
there exists a measurable set such that the behavior of its
characteristic function, regarded as a vector of length
$2^n$ with coordinates $0$, $1$, in the hierarchy of conditional
measures corresponding to the first $n$ partitions does not stabilize
even for the exponential scaling, contradicting the standardness
of the filtration. The calculations carried out in the paper
not only proved that the filtration is nonstandard, but also
gave a lower bound on the exponential entropy in this example.
This was the first application of the standardness criterion
and a motivation for introducing exponential entropy.
The same bound applied to the exponential entropy of the
dyadic filtrations arising from Bernoulli actions of infinite
sums of the groups $\Bbb Z_k$, implying the nonstandardness
of these filtrations.
In the latter case, the exponential entropy numerically coincides
with the entropy of the action, though their definitions
are quite different (see \cite{V5, He1}).

In the class of $r$-adic filtrations, natural examples of stationary
filtrations arise from random walks over the trajectories of automorphisms,
with equiprobable transitions to one of the $r$
points of the trajectory; for example,
the so-called $(T,T^{-1})$-endomorphism is
the random walk that moves from a point
$x$ to the points $Tx$ and $T^{-1}x$ with probabilities
$1/2$; in the more general case of an $r$-adic filtration,
the random walk moves with probability $1/r$ to one of the $r$
points of the trajectory of an arbitrary group action. In general,
random walks (and, more generally, the theory of polymorphisms)
provide many interesting examples of filtrations. The first example
given above also belongs to this class.

Much later, another example of a random walk over
trajectories was given, also using the standardness criterion, this time
in the positive direction: the filtration of pasts
of the $(T,T^{-1})$ random walk constructed from the rotation
$T_{\lambda}$ of the circle by an irrational angle
$\lambda$ is standard; this was established first for
values of $\lambda$ that can be well approximated by rational numbers
\cite{V6, Parry}, and then for arbitrary values of
$\lambda$ \cite{HR}.

The conjecture of the first author
that the filtration of pasts of the random walk over the trajectories
of a Bernoulli action of $\Bbb Z$ (the so-called Kalikow endomorphism)
is nonstandard
had been open for a long time. S.~Kalikow  \cite{K} showed that the
$(T,T^{-1})$ endomorphism is not even
loosely Bernoulli. This was the first example of a
natural non-Bernoulli endomorphism.
The question naturally arose about the type of the
filtration of pasts of this endomorphism.
It is not difficult to check that its exponential entropy
vanishes. Finally, in
\cite{HH} it was proved, with the help of the standardness criterion, that
the filtration of pasts of this endomorphism is indeed nonstandard.
In fact, the proof explicitly used scaled entropy, see below.

In \cite{H1,Fel}, examples are constructed showing that the Bernoulli
property of an endomorphism and the standardness of the filtration of pasts
are in general position.
In other words,

1) Bernoulli automorphisms have generators leading to random processes
with nonstandard filtrations of pasts; such are, for example, the above
random walk over the trajectories of a Bernoulli action of the free group,
and the examples of random walks considered below;

and

2) there exist stationary random processes with standard filtrations of pasts
such that the shifts in the spaces of trajectories of these processes
are not isomorphic to a Bernoulli shift
(\cite{Fel}).

\subsection{Scaling for random walks over the trajectories of Bernoulli
group actions}

The next class of examples of filtrations is generated by random walks
over the trajectories of Bernoulli actions of arbitrary groups.
Assume that we are given an arbitrary countable group
$G$ with finitely many generators
$g_1, \dots, g_s$, and the Bernoulli action of this group by left shifts
$T_{g_i}$ in the space $F$ of all $\{0,1\}$-functions on $G$
endowed with the product measure with the factor
$(1/2,1/2)$ (i.e., a Bernoulli measure). In what follows,
the set of possible values
of functions in $F$ and its (finite) cardinality
are irrelevant, so that for simplicity we
restrict ourselves to the values
$0$ and $1$. Consider the random walk on the space $F$ over the trajectories of
the action of the group $G$ with equal transition probabilities:
$$
{\rm Prob}(f(x)\mapsto f(g_i^{\pm 1}x))={1\over 2s}.
$$

Thus we consider a generalization of the
$(T,T^{-1})$ construction in which the group
$\Bbb Z$ is replaced with an arbitrary discrete group $G$
with finitely many generators. In the previous notation, $X$ is
the space of trajectories of the Markov process with the state space
$F(G,\{0;1\})=2^G$, the Bernoulli measure on
$F$ as an invariant measure, and the above transition probabilities. This Markov process will
be called the {\it random walk over the trajectories of the Bernoulli action} of $G$.
The sequence of pasts of this Markov process is an $r$-adic filtration with
$r=2s$. What can we say about this filtration, in particular,
about its scaled entropy and scaling? In full generality this problem
is far from being solved; note that the first example of a nonstandard
filtration was exactly of this kind, with the free group as $G$.
Below we give its solution for lattices and nilpotent groups.

The paper \cite{RHH} in fact provides a bound on the scaling function
(though the authors of \cite{RHH} do not use the entropy terminology)
for the filtration of pasts of Kalikow's
$(T,T^{-1})$ endomorphism, where $T$ is a Bernoulli automorphism. As shown in
\cite{G2}, the scaling function in this case is equivalent to
$c(\varepsilon,n)=(n\log({1\over\varepsilon}))^{1/2}$.

The paper \cite{STHO} generalized the problem about the
$(T,T^{-1})$ automorphism
solved by S.~Kalikow. Namely, the authors of
that paper considered the Markov automorphism of a simple walk
over the trajectories of an action of the lattice
${\Bbb Z}^d$. But they were interested not in the
filtration, but in the type of the Markov shift. They showed
that the situation depends crucially on the dimension $d$ of the lattice:
while in \cite{Kal} (for $d=1$) it was proved that the
$(T,T^{-1})$ automorphism is not even
loosely Bernoulli and, all the more so, not Bernoulli, for
$d>1$ it is Bernoulli (though the quality of the natural generator,
the type of Bernoulli property, depends on $d$).

We study the scaled entropy and, in particular, the scaling
of the corresponding filtration for the groups
${\Bbb Z}^d$; in this case, it turns out that the entropy
properties depend on the dimension of the lattice
in a quite regular way.

\begin{theorem}[see \cite{G4}] For the group $G={\Bbb Z}^d$
the proper scaling function of the filtration is of
polynomial growth:
$c(\varepsilon,n)=(n\log({1\over\varepsilon}))^{d/2}$.
\end{theorem}

A further generalization of combinatorial techniques allowed the second
author to find the proper scaling function for the case of random
walks over the trajectories of a Bernoulli action of an arbitrary
countable nilpotent group $G$. Recall that the weighted rank
(in the continuous case, the weighted, or Hausdorff, dimension)
of a nilpotent group is the number
$d=\sum i(n_i-n_{i-1})$, where $(n_1,{\ldots} ,n_k)$ is the vector of ranks
of the groups $H_i$, with $H_i$ being the quotient
by the $i$th element of the lower central series of the group $G$.

\begin{theorem}
For a countable nilpotent group $G$ and the filtration
$\Xi=\{\xi_n\}_{n=1}^{\infty}$ of pasts of the Markov random walk
over the trajectories of a Bernoulli action of $G$, the proper scaling function
is equivalent to
$$
c(\varepsilon,n)=(n\log({1\over\varepsilon}))^{d/2},
$$
where $d=d(G)$ is the weighted rank of $G$.
\end{theorem}

Thus for abelian and nilpotent groups the scaling is polynomial;
for a walk over the trajectories of a Bernoulli action of a free nonabelian
group, the scaling is exponential, and the entropy of the
corresponding filtration is the ordinary (exponential) entropy.
The question about the scaling for random walks on solvable groups is open.

\subsection{Necessary estimates}

Here we will give only several separate statements that constitute
the main part of the proof of the theorems on scaling, leaving the
detailed exposition till another occasion. The first result needed
for estimating the entropy of the iterated metric generalizes a result from
\cite{HH} and can be proved in a similar way.

\begin{theorem} There exists a semimetric $\rho$ on the space $X$ and
a number $\varepsilon>0$ such that for every polynomial $p$ there exists
$n_0$ such that for all points $x$ except a set of measure
$\varepsilon$,
$$
\mu\{y : \rho_n(x,y)<\varepsilon\}<{1\over p(n)},\qquad n>n_0,
$$
where $\rho_n$ is the iterated semimetric (see \S\,2).
\end{theorem}

The second result is a generalization of a technical result from
\cite{RHH,Kal}.

\begin{theorem} Let $X$ be the space of trajectories of the Markov process
of a random walk over the trajectories of a Bernoulli action (see above)
of a nilpotent group $G$ with generators
$e_i$, $i=1, \dots, n$. For every
$\delta>0$ we can find a subset $M \subset
X$ of measure $1-\delta$ and a number $h_0$ such that for every
$h>h_0$ and every pair of trajectories $\{u_i\},\{v_i\}$ from $M$,
there exists $n\in [h,h^5]$ such that
$$
{1\over \sqrt{n}}\|\prod\limits_{i=1}^{n} u_i\|<c\quad\hbox{  and  }\quad
{1\over \sqrt{n}}\|\prod\limits_{i=1}^{n} u_i\|<c,
$$
where $\|\cdot \|$ is understood as the number of factors in the minimal
representation of a group element as a product of the generators
$g_i$ and their inverses.
\end{theorem}

Note that the set $M$ determined by the conditions of the previous theorem
can be chosen in different ways, but in what follows only the existence
of such a set is of importance.

It will be convenient to slightly modify the construction of the space $X$
and represent the shift in this space as a skew product.
Namely, the new version of $X$ is
$X=F(G)\times
B^{\infty}$, and the Markov shift $T$ in $X$ is a skew product
over the one-sided Bernoulli shift.

Let us choose an initial semimetric
$\rho$ on the space $X=F(G)\times B^{\infty}$ that is measurable
with respect to the partition $\eta_m$ into cylinders of order $m$
in the sense of the structures of the spaces  $F(G)$ and
$B^{\infty}$. Note that with this choice of a semimetric,
the corresponding iterated semimetric
$\rho_n$ will be measurable with respect to  $\eta_{m+n}$;
indeed, the preimage $T^{-1}(\eta_m)$ of $\eta_m$ is measurable with respect to
$\eta_{m+1}$.

In order to estimate the scaling of the filtration of pasts of the Markov
process, we estimate the
$\ve$-entropy of the measure in the space
$(X,\mu,\rho_n)$. It will be convenient to use the combinatorial
description of the iterated semimetric similar to that given in Sec.~2.2.

The iterated semimetric $\rho_n$ can be expressed in terms of the
initial metric and the group of automorphisms of the tree as follows:
$$
\rho_n(x,y)=\min_{a\in D_{n,r}}{1\over r^n}\sum\limits_{i=1}^{r^n}
\rho(x_i,y_{a(i)}).
$$

This explicit expression allows us to estimate the values of the iterated
metric using the properties of the group
$D_{n,r}$. For example, to estimate the entropy from above, it suffices
to estimate the iterated metric from above. To this end, simply
replace the group
$D_{n,r}$ by the trivial group, immediately obtaining the bound
$$
\rho_n(x,y)\le {1\over r^n}\sum\limits_{i=1}^{r^n} \rho(x_i,y_i),
$$
which, after applying the central limit theorem for the nilpotent group
$G$ (see \cite{Raugi}) gives the required upper bound.

The lower bound on the entropy (and, correspondingly, on the iterated
semimetric) is slightly more difficult to obtain. The iterated semimetric
does not allow for a simple uniform bound, but we can estimate the measure
of a typical $\ve$-ball in the iterated metric
$\rho_n$. Such a bound also allows us to estimate the
$\ve$-entropy.

The fact that the semimetric
$\rho_n$ is measurable with respect to the cylinder partition allows us
to translate the problem into purely combinatorial terms:
we need to estimate the number of cylinders of order
$n+m$ lying at a distance less than
$\ve$ from a given cylinder on the quotient space
$X/\eta_{n+m}$ with the reduced semimetric
(see Sect.~2.2).  Recall that $D_{n,r}\subset S_{r^n}$, where $S_{r^n}$
is the symmetric group acting by substitutions on the space
$K_n=\{0,1\}^{2^n}$ of sequences of $0$'s and $1$'s of length $n$.
Let us fix two typical points (i.e., points from the set $M$)
$x$ and $y$ lying at a distance at least
$\ve$. The proof reduces to studying the properties of the automorphism
$a\in D_{n,r}$ that realizes the distance between these points. Consider
the projection  $x \mapsto \bar x$ of the space $X$ to the quotient space
$X/\eta_{n+m}$ and the embedding
$* : \bar{x}\mapsto \bar x^*=\{x_i,i=\overline{1,{\ldots} ,r^n}\}$
that sends a point to the collection of its preimages. The main property
of the automorphism $a$ is that it superposes two vectors,
$\bar{x}^*\in K_n$ and
$\bar{y}^* \in K_n$, of exponential (in $n$) length
$r^n$, and these vectors are in turn determined by vectors,
$\bar{x},\bar{y}\in
X/\eta_{n+m}$, of polynomial length
$(n+m)^r$.

It turns out that the action of the minimizing automorphism $a$
``almost'' factorizes. This means that the diagram
$$
\begin{array}{ccc}{K_n}& \stackrel{a}{\longrightarrow}& {K_n}\\
{\Big\uparrow}* & &
{\Big\uparrow} * \\X/\eta_{n+m}&\stackrel{A}{\longrightarrow}& X/\eta_{n+m}\\
\end{array}
$$
is almost commutative, in the sense that the (Hamming) distance
between the results of following two paths in this diagram
does not exceed some $\delta$ which tends to zero as $n$ tends
to infinity:
$$d(a(\bar{x}^*), A(\bar{x})^*)<\delta.$$

In addition to the fact that the automorphism admits the above
``$\delta$-factorization,'' it turns out that the constructed quotient map
will be almost identical. Further calculations show that the number of classes
$\bar{y}\in X/\eta_{n+m}$ close to $\bar{x}$ in the iterated semimetric $\rho_n$
is subexponential, i.e., inessential from the point of view of the
$\ve$-entropy.

\smallskip
Translated by N.~V.~Tsilevich.


\begin{thebibliography}{}
\bibitem{V0} A.~M.~Vershik, Decreasing sequences of measurable partitions and their applications,
{\em Sov. Math. Dokl.} {\bf 11} (1970),
1007--1011.

\bibitem{V3}A.~M.~Vershik, Continuum of pairwise
nonisomorphic dyadic sequences,
{\em Funct. Anal. Appl.}  {\bf 5}
(1971), 182--184.

\bibitem{V2}A.~M.~Vershik, Approximation in measure theory,
D.Sc. Thesis, Leningrad, 1973.

\bibitem{V5}A.~M.~Vershik, Four definitions of the scale of an automorphism,
{\em Funct. Anal. Appl.} {\bf 7} (1973), 169--181.

\bibitem{V1}A.~M.~Vershik, Theory of decreasing sequences of
measurable partitions,
{\em   St.~Petersburg Math. J.} \textbf{6}, No. 4 (1995), 705--761.

\bibitem{V4}A.~M.~Vershik,
Dynamic theory of growth in groups: entropy, boundaries, examples,
{\em Russian Math. Surveys} {\bf 55}, No. 4
(2000), 677--733.

\bibitem{V7}A.~M.~Vershik, The Kantorovich metric: the initial history and
little-known applications,
{\em  J. Math. Sci. (New York)} {\bf 133}, No. 4 (2006), 1410--1417.

\bibitem{G1} A.~D.~Gorbulsky, On an entropy property of a decreasing
sequence of measurable partitions, {\em  J. Math. Sci. (New York)}
{\bf 107}, No. 5 (2001), 4157--4160.

\bibitem{G2} A.~D.~Gorbulsky, A correlation between different definitions of the entropy of decreasing
sequences of partitions; scaling,
{\em  J. Math. Sci. (New York)}
{\bf 121}, No. 3 (2004), 2319--2325.

\bibitem{G3} A.~D.~Gorbulsky, An example of the interpolation of the standard sequence of measurable
partitions with large scaling entropy,
{\em  J. Math. Sci. (New York)}
{\bf 126}, No. 2 (2005), 1043--1045.

\bibitem{G4} A.~D.~Gorbulsky,
The $\sigma$-algebra of the pasts of a random walk on the orbits of
the Bernoulli action of the group $Z^D$,
{\em  J. Math. Sci. (New York)}
{\bf 138}, No. 3 (2006), 5686--5690.

\bibitem{K} L.~V.~Kantorovich, On the translocation of masses,
{\em Dokl. Akad. Nauk SSSR} {\bf 37}, Nos.~7--8 (1942), 227--229.

\bibitem{Ku} A.~G.~Kushnirenko, On metric invariants of entropy type,
{\em Russian Math. Surveys} {\bf22}, No. 5 (1967), 53--61.

\bibitem{R1} V.~A.~Rokhlin, On the fundamental ideas of measure theory,
{\em Mat. Sbornik N. S.} {\bf 25(67)} (1949), 107--150.

\bibitem{R2} V.~A.~Rokhlin, Lectures on the entropy theory of
measure-preserving transformations,
{\em Russian Math. Surveys} {\bf22}, No. 5 (1967), 1--52.

\bibitem{S}A.~M.~Stepin, On entropy invariants of decreasing sequences of measurable partitions,
{\em Funct. Anal. Appl.} {\bf 5}
(1971), 237--240.

\bibitem{BR} X.~Bressaud, A.~Maass, S.~Martinez, and J.~San Martin,
Stationary processes whose filtrations are
standard, {\em Ann. Probab.} {\bf 34},
No. 4 (2006), 1580--1600.

\bibitem{Raugi} P.~Cr\'epel and A.~Raugi,
Th\'eor\`eme central limite sur les groupes nilpotents,
{\em Ann. Inst. H. Poincar\'e Sect. B (N. S.)} \textbf{14}
(1978), 145--164.

\bibitem{DFST} L.~Dubins, J.~Feldman, M.~Smorodinsky,
and B.~S.~Tsirelson, Decreasing sequences of $\sigma$-fields and a measure
change for Brownian motion, {\em Ann. Probab.} {\bf 24}, No. 2
(1996),  905--911.

\bibitem{Em} M.~Emery, Espaces probabilis\'es filtr\'es: de la th\'eorie de
Vershik au mouvement brownien, via des id\'ees de Tsirelson, {\em
Seminaire Bourbaki}, Vol. 2000/2001, {\em Ast\'erisque},
No. 282  (2002), Exp. No. 882, vii, 63--83.

\bibitem{Fel} J.~Feldman and D.~Rudolph, Standardness of the decreasing sequences of $\sigma$-fields given by certain dyadic endomorphisms,
 {\em Fund. Math.} {\bf 157}, Nos. 2--3 (1998), 175--189.

\bibitem{He1} D.~Heicklen, Entropy and $r$ equivalence,
{\em Ergodic Theory Dynam. Systems} {\bf  18}, No. 5  (1998), 1139--1157.

\bibitem{HH} D.~Heicklen and  C.~Hoffman, $T,T^{-1}$ is not
standard, {\em Ergodic Theory Dynam. Systems} {\bf 18}, No. 4 (1998),
875--878.

\bibitem{RHH} D.~Heicklen, C.~Hoffman, and D.~Rudolph, Entropy and
dyadic equivalence of random walks on a random scenery, {\em
Adv. Math.} {\bf 156}, No. 2 (2000), 157--179.

\bibitem{H1} C.~Hoffman and D.~Rudolph, A dyadic endomorphism which is Bernoulli but not standard,
{\em Israel J. Math.}  {\bf130} (2002), 365--379.

\bibitem{HR} C.~Hoffman and D.~Rudolph, Uniform endomorphisms which are
isomorphic to a Bernoulli shift, {\em Ann. Math.} (2),  {\bf 156}, No. 1
(2002),  79--101.

\bibitem{STHO}  F.~Hollander and J.~Steif, Random walk in random
scenery, {\em IMS Lect. Notes} {\bf 48} (2006), 53--65.

\bibitem{Kal} S.~A.~Kalikow, $T,T^{-1}$ transformation is not loosely
Bernoulli, {\em Ann. Math.} {\bf 115} (1982), 393--409.

\bibitem{Parry} W.~Parry, Automorphisms of the Bernoulli endomorphism and a class of skew-products.
{\em Ergodic Theory Dynam. Systems} \textbf{16}, No. 3 (1996),  519--529.

\bibitem{Sh} P.~Shields and K.~Marton, How many future measures can there be?,
 {\em Ergodic Theory Dynam. Systems} {\bf 22} (2002), 257--280.

\bibitem{V6}A.~M.~Vershik, Pasts of $T,T^{-1}$ are nonstandard,
manuscript. Berkeley, 1995.

\end{thebibliography}
\end{document}